\documentstyle[12pt]{article}

  \newcommand{\const}{\rm const}

  \newcommand{\sign}{\rm sign}

\begin{document}

   \begin{center}

{\bf Modulus of continuity for superlacunar trigonometric series  and} \par

\vspace{3mm}

{\bf continuity of Gaussian  stationary random processes.  }\par

\vspace{4mm}

{\bf M.R.Formica,  E.Ostrovsky, and L.Sirota.}

  \end{center}

\vspace{4mm}

 Universit\`{a} degli Studi di Napoli Parthenope, via Generale Parisi 13, Palazzo Pacanowsky, 80132,
Napoli, Italy. \\

e-mail: mara.formica@uniparthenope.it \\

\vspace{4mm}

 \ Department of Mathematics and Statistics, Bar-Ilan University, \\
59200, Ramat Gan, Israel. \\

e-mail: eugostrovsky@list.ru\\

\vspace{4mm}

\ Department of Mathematics and Statistics, Bar-Ilan University, \\
59200, Ramat Gan, Israel. \\

e-mail: sirota3@bezeqint.net \\

\vspace{4mm}

\begin{center}

 \ {\bf Abstract}

\vspace{5mm}

 \end{center}

 \  We deduce bilateral interrelations between Fourier coefficients for lacunar trigonometric series and
 modulus of their continuity. We obtain also  as an application some conditions for continuity and discontinuity for
 Gaussian periodic stationary random centered processes. \par

\vspace{5mm}

 \ {\it Key words and phrases.} Lacunar and superlacunar trigonometric series and functions, necessary condition and 
sufficient one,  modulus of continuity,  Fourier coefficients, ordinary and complex functions, Gaussian random variable 
and stationary random  processes, convergence, estimation. \par

\vspace{5mm}

\section{Introduction. Lacunar and superlacunar trigonometric series.}

\vspace{5mm}

 \hspace{3mm} Let $ \  \vec{n}  = \{n(k)\}, \ k = \pm 1, \pm 2, \pm 3,\ \ldots  \  $ be a set of integer numbers, positive as well
as negative,  such that  $ \ n(1) \ge  0, \ n(-1) \le 0, \ $

\begin{equation} \label{set nk}
n(k+1) \ge n(k)+1, \ k \ge 1; \hspace{3mm}  n(k-1) \le n(k) - 1, \ k < 0.
\end{equation}

 \ Let also $ \   \vec{c} = \{c(k) \},   \ k \in \vec{n}  \  $ be a numerical  absolute summable numerical sequence

\begin{equation} \label{abs converg}
\sum_{k \in \vec{n}} |c(k)| < \infty.
\end{equation}

 \ Define the following  $ \ 2 \pi \ $   periodical continuous function  by means of  absolute convergent Fourier series

\begin{equation} \label{key series}
f(t)  = f[\vec{n}, \vec{c}](t) \stackrel{def}{=} \sum_{k \in \vec{n}} c(k) \ \exp( i \ t \ n(k)). \ t \in [-\pi, \pi].
\end{equation}

 \ Of course,

$$
c(k)  = (2 \pi)^{-1} \int_{-\pi}^{\pi} \ \exp( - i \ t \ n(k)) \ f(t) \ dt.
$$

\vspace{3mm}

 \  Recall that the series of the form   (\ref{key series}) are named {\it lacunar,} iff

\vspace{3mm}

\begin{equation} \label{lacunar ser plus}
 \underline{\lim}_{ k \to \infty}  \left[ \ \frac{n(k + 1)}{n(k)} \ \right] > 1.
\end{equation}

and

\begin{equation} \label{lacunar ser minus}
 \underline{\lim}_{ k \to \infty}  \left[ \ \frac{|n(-k - 1)|}{|n(-k)|} \ \right] > 1.
\end{equation}

\vspace{4mm}

 \ {\bf Definition 1.1.}   The series of the form   (\ref{key series}) are named {\it superlacunar,} iff

\vspace{3mm}

\begin{equation} \label{superlacunar ser plus}
\underline{\lim}_{ k \to \infty} \left[ \ \frac{n(k + 1)}{n(k)}  \ \right] = \infty
\end{equation}

and  in addition

\begin{equation} \label{superlacunar ser minus}
\underline{\lim}_{ k \to \infty} \left[ \ \frac{|n(-k - 1)|}{|n(-k)|}  \ \right] = \infty,
\end{equation}
see e.g.  \cite{Bary} , chapter 5; \cite{Zygmund}, chapter 3. \par

\vspace{4mm}

 \hspace{3mm} {\bf Our claim in this report is to establish the interrelations between modulus of continuity of the
 superlacunar function and  behavior of  its Fourier coefficients. } \par

\vspace{3mm}

 \ {\bf We consider as an application the continuity properties of Gaussian periodic stationary random processes. } \par

 \ As an example: the famous Weierstrass function. \par

\vspace{5mm}

\section{Bilateral inequalities. Examples.}

\vspace{5mm}

 \hspace{3mm} Recall that the modulus of continuity $ \ \omega[f](\delta), \ \delta \in [0, 2\pi] \ $  for the continuous $ \ 2 \pi \ $
periodical numerical valued  function $ \ f = f(t) \ $  is defined as follows

\begin{equation} \label{modulus cont}
\omega[f](\delta) \stackrel{def}{=}  \sup_t \ \sup_{h: |h| \le \delta} | \  f(t + h) - f(t)  \ |.
\end{equation}

  \ It is well - known  \ \cite{Zygmund}, chapters 2,3 \ that for the function in particular of the form  (\ref{key series})

\begin{equation} \label{coeff estim}
|c(k)| \le (4 \pi)^{-1} \ \omega[f](\pi/n(|k|)), \  |k| \ge 1.
\end{equation}

\vspace{4mm}

 \ Conversely, let the function $ \ f = f(t), \ t \in R \ $ has the form  (\ref{key series}), such that (we recall)
 $ \ \vec{c} \in l_1. \ $ Let also $ \ t,s \in T, \ |t - s| < h, \ h \in [0, 2 \pi]. \ $  We have

$$
|f(t) - f(s)| \le h \ ||\vec{c}||l_1 \  \sum_{k: |k| \le N} |k| \ |n(|k|)|  +
$$

$$
 2 \ \sum _{k: |k| > N} |c(k)| \ = h \Sigma_1[c,n](N) + \Sigma_2[c](N),
$$
where

$$
||\vec{c}||l_1 = \sum_{k} |c(k)| < \infty, \ \Sigma_1[c,n](N) = ||\vec{c}||l_1 \cdot  \sum_{k: |k| \le N} |k| \ |n(|k|)|,
$$

$$
\Sigma_2[c](N) = 2 \ \sum _{k: |k| > N} |c(k)|; \ N \ge 2.
$$

\vspace{3mm}

 \ To summarize  common with (\ref{coeff estim}):  \par

\vspace{4mm}

 \ {\bf Proposition 2.1.} Denote

\vspace{3mm}

$$
\Sigma[c, \vec{n}] (\delta) \stackrel{def}{=} \inf_{N = 2,3,4,\ldots } \left[ \ \delta \cdot \Sigma_1[c,n](N) + \Sigma_2[c](N) \ \right].
$$

 \ We get under formulated above  notations and conditions

\vspace{4mm}

\begin{equation} \label{key est}
\omega[f](\delta) \le \Sigma[c, \vec{n}] (\delta), \ \delta \in [0, 2 \pi).
\end{equation}

\vspace{4mm}

 \ {\bf Remark 2.1.} Obviously,

$$
\lim_{\delta \to 0+} \Sigma[c, \vec{n}] (\delta) = 0.
$$

\vspace{3mm}

\ Let us consider some examples. \par

\vspace{3mm}

 \ {\bf Example 2.1.} Let us put

$$
  c(k) = k^{-2 \Delta}, \ k \ge 1;  \ c(k) = 0, \ k \le 0;  \ n(k) \stackrel{def}{=} 2^k, \ \Delta = \const > 1/2,
$$
the lacunar case.  \ On the other words, we consider  the following complex valued continuous function

\begin{equation} \label{first fun}
g(t) = g_{\Delta}(t) := \sum_{k=1}^{\infty} k^{-2\Delta} \ \exp \left( \  2^k \ i \ t \ \right).
\end{equation}

\vspace{3mm}

 \ It follows from (\ref{coeff estim})

$$
\omega[g](\delta) \ge C_1 \ |\ln \delta|^{- 2 \Delta}, \ 0 < \delta  < 1/e, \ 0 < C_1 < \infty,
$$
a lower bound. The correspondent upper bound may be found on the basis of proposition 2.1:

$$
\omega[g](\delta) \le C_2 \ |\ln \delta|^{1 - 2 \Delta}, \ 0 < \delta  < 1/e, \ 0 < C_2 < \infty.
$$

\vspace{4mm}

\ {\bf Example 2.2.} Let us introduce now

$$
  c(k) = k^{-2 \Delta}, \ k \ge 1;  \ c(k) = 0, \ k \le 0;  \ n(k) \stackrel{def}{=} 2^{2^k}, \ \Delta = \const > 1/2,
$$
the superlacunar case.  \  We consider   the following complex valued continuous function

\begin{equation} \label{second fun}
s(t) = s[\Delta](t) := \sum_{k=1}^{\infty} k^{-2\Delta} \ \exp \left( \ i \cdot 2^{2^k} \cdot t \ \right).
\end{equation}

\vspace{3mm}

 \ It follows from (\ref{coeff estim})  as before

$$
\omega[s](\delta) \ge C_3 \ \left[ \ \ln |\ln \delta|  \ \right]^{- 2 \Delta}, \ 0 < \delta  < e^{-e}, \ 0 < C_3 < \infty,
$$
a lower bound. The correspondent upper bound may be found  alike on the basis of proposition 2.1:

$$
\omega[s](\delta) \le C_4 \  \left[ \ \ln |\ln \delta| \ \right]^{1 - 2 \Delta}, \ 0 < \delta  < e^{-e}, \ 0 < C_4 < \infty.
$$

\vspace{5mm}

\section{Gaussian periodic stationary  random processes.}

\vspace{5mm}

 \hspace{3mm}  Let $  \ \eta(t), \ t \in  T := [ - \pi,\pi] \ $ be centered (mean zero) stochastic continuous  real valued
 separable Gaussian  distributed random process  (r.p.) having continuous covariation function

$$
R(t,s)  := {\bf E} \eta(t) \eta(s), \ t,s \in T.
$$

\vspace{3mm}

 \ The problem of finding   the conditions (necessary conditions and sufficient ones) for the continuity with probability 
 one of the r.p.  $ \ \eta(\cdot): \ $

\begin{equation} \label{statement}
{\bf P}(\eta(\cdot) \in C(T)) = 1
\end{equation}
can be considered  as a classic, see e.g. \cite{Belyaev},  \cite{Buldygin}, \cite{Dudley}, \cite{Fernique}, \cite{Garsia},
\cite{Kozachenko-Ostrovsky 1985}, \cite{Leadbetter}, \cite{Marcus 1}  - \cite{Marcus 6}, \cite{Ostrovsky 0} etc. We mention
here  the famous result belonging to X.Fernique  \ \cite{Fernique}: if the following integral convergent

\begin{equation} \label{cond Fernique}
\int_0^{\infty} \omega^{1/2}[R](\exp(-x^2/2)) \ dx < \infty,
\end{equation}
then the relation   (\ref{statement}) holds true. So, the condition (\ref{cond Fernique}) is {\it sufficient} for the continuity
a.e. of the Gaussian r.p. $ \ \eta(t). \ $  Open question: what happens if the condition  (\ref{cond Fernique}) is not satisfied? \par

\vspace{4mm}

 {\bf  We will show in this report that  if the condition  (\ref{cond Fernique}) is not satisfied, i.e.  when the covariation function
 of the Gaussian  r.p. may be  arbitrarily non - smooth (as rough as you like), then the  correspondent  stationary Gaussian r.p.
may be continuous or conversely may be extremely discontinuous, i.e. may be unbounded on every non - empty interval.} \par

\vspace{4mm}

 \ {\bf We consider only the case of stationary and periodical centered Gaussian r.p. } \par

\vspace{4mm}

 \ To be more concrete, introduce the following Gaussian distributed stationary periodical process of the form  (super - lacunar series)

\begin{equation} \label{key rand process}
\nu(t) = \nu_{\Delta}(t) \stackrel{def}{=} \sum_{|k| \ge 1} |k|^{-\Delta} \ \kappa_k \ \exp \left( \ i \cdot \sign(k) \cdot  2^{2^|k|}  \ \right),
\end{equation}

$$
\Delta = \const > 0, \ t \in R,
$$
and $ \ \{\kappa_k\}  \ $  are independent  centered real valued  standard Gaussian distributed  random variables: $ \  {\bf E} \kappa^2_k = 1.  \ $ \par

 \vspace{3mm}

  \ The covariation function $ \  r(t) = {\bf E} \xi(t+s) \ \overline{\xi}(s)  \ $   of the stationary r.p. $ \ \nu(t) \ $  has a form (superlacunar series)

\vspace{3mm}

\begin{equation} \label{covar fun}
r(t) = 2 \sum_{k=1}^{\infty} k^{-2 \Delta} \ \cos \left( \  2^{2^k} \ t  \ \right), \ t \in R.
\end{equation}

\vspace{3mm}

 \ If $ \  \Delta \le 1/2, \ $  then the function $ \ r = r(t) \ $ is  extremely discontinuous.
 \ If $ \ \Delta \le 1,  \ $ then the r.p. $ \ \nu = \nu(t) \ $ is  also extremely discontinuous a.e., \ \cite{Belyaev}. \par

 \vspace{3mm}

  \ The behavior of the function $ \ r = r(t) \ $ in the case  when $ \ \Delta > 1 \ $  is obtained in the example 2.2.
 So, it is very  rough, despite the  correspondent r.p.  $ \ \nu(t) \ $ is continuous. \par
  \ Moreover, the covariation function $ \ r(t) \ $  for {\it continuous} Gaussian stationary  centered random process may be
 as rough as you like, for instance  when

 \vspace{3mm}

\begin{equation} \label{covar fun  very}
r(t) = 2 \sum_{k=1}^{\infty} k^{-2 \Delta} \ \cos \left( \ n(k) \ t  \ \right), \ t \in R, \ \Delta = \const > 1,
\end{equation}
where as before $ \  \{n(k)\} \ $ is very fast increasing integer sequence, see (\ref{coeff estim}). \par

\vspace{5mm}

\section{Concluding remarks.}

\vspace{5mm}

 \hspace{3mm} An open question: let $ \  \omega = \omega(\delta), \ \delta \in [0,1/e]   \ $ be arbitrary non - trivial modulus of continuity.
 There exists or not a stationary Gaussian random process  $ \ \zeta = \zeta(t), \ t \in R  \ $ with continuous sample path
for which its covariation function $ \ r =  r[\zeta](t)  \ $ obeys a property

$$
\omega[r](\delta) \asymp \omega(\delta), \ \delta \in [0,1/e].
$$

\vspace{6mm}

\vspace{0.5cm} \emph{Acknowledgement.} {\footnotesize The first
author has been partially supported by the Gruppo Nazionale per
l'Analisi Matematica, la Probabilit\`a e le loro Applicazioni
(GNAMPA) of the Istituto Nazionale di Alta Matematica (INdAM) and by
Universit\`a degli Studi di Napoli Parthenope through the project
\lq\lq sostegno alla Ricerca individuale\rq\rq .\par

\end{document}